\documentclass[12pt]{article} 
\usepackage{amssymb,ifthen,graphicx}
\renewcommand{\thefootnote}{\fnsymbol}
\newcommand\blfootnote[1]{%
  \begingroup
  \renewcommand\thefootnote{}\footnote{#1}%
  \addtocounter{footnote}{-1}%
  \endgroup
}

\renewcommand{\phi}{\varphi}
\newcounter{rcc}
\newcommand{\rc}{\roman{rcc}\stepcounter{rcc}}

\newcommand{\be}{\begin{eqnarray*}}
\newcommand{\ee}{\end{eqnarray*}}
\newcommand{\pf}{\noindent {\em (Proof)}}
\newcommand{\qed}{\hfill\fbox{\rule[-1mm]{0mm}{1mm}\,}}

\newcommand{\wt}[1]{\widetilde{#1}}
\newcommand{\maps}[2]{\! :\! #1\rightarrow\! #2}
\newcommand{\g}[2]{\langle#1,#2\rangle_g}

\newcommand{\thmref}[1]{Theorem~\ref{#1}}

\newcommand{\lemref}[1]{Lemma~\ref{#1}}

\newtheorem{lemma}{Lemma}[section]

\newtheorem{theorem}[lemma]{Theorem}

\newtheorem{defn}[lemma]{Definition}

\begin{document}

\title{Hyperbolic metrics on open
subsets of Ptolemaic spaces with sharp parameter bounds}

\author{Neil N.\,Katz}

\date{}
\maketitle

\begin{abstract}
It is shown that a construction of Z.\,Zhang and Y.\,Xiao on open subsets
of Ptolemaic spaces yields, when the subset has boundary containing at least
two points, metrics that are Gromov hyperbolic with parameter $\log 2$ and
strongly hyperbolic with parameter $1$ with no further conditions on the
open set.  A class of examples is constructed on Hadamard manifolds showing
these estimates of the parameters are sharp.
\end{abstract}

\section*{Introduction}

\blfootnote{2010 Mathematics Subject Classification: 51M10, 53C23.}
\blfootnote{keywords: Ptolemaic, Gromov hyperbolic, strongly hyperbolic,
metric space.}

In this paper a construction in \cite{ZX} is applied to produce a metric on
an open subset of a Ptolemaic space.  When the open set has boundary
containing at least two points, the metric is strongly hyperbolic with
parameter 1 and therefore (by Theorem 4.2 in \cite{NS}) it is Gromov
hyperbolic with parameter $\log 2$.  This is shown with no other conditions
on the boundary (\thmref{maintheorem}).  A class of examples is constructed
for an open subset of a Hadamard manifold.  (By Theorem 1.1 in \cite{BFW} a
complete Riemannian manifold is Ptolemaic if and only if it is Hadamard.) In
the examples the open set has a boundary with at least two isolated points,
one of which is the nearest neighbor in the boundary to the other
(\thmref{hadamard}).  This class of examples are all Gromov hyperbolic with
parameter $\log 2$ and no lower, so by applying the same theorem \cite{NS}
and \thmref{maintheorem} they are strongly hyperbolic with parameter 1 and
no higher.

\begin{defn}
Let $(X,d)$ be a metric space.
  
It is Ptolemaic iff
$$
d(x,y)\,d(z,t)\leq d(x,z)\,d(y,t)+d(x,t)\,d(y,z)
$$
for any $x,y,z,t\in X$.

It is hyperbolic in the sense of Gromov {\rm (\cite{G}; see also \cite{JV})}
with parameter $\delta>0$ iff
$$
d(x,y)+d(z,t)\leq \max\left\{d(x,z)+d(y,t),d(x,t)+d(y,z)\right\}+2\delta
$$
for any $x,y,z,t\in X$.

It is strongly hyperbolic {\rm (\cite{NS})} with parameter $\epsilon>0$ iff
\be
& & \exp\left(\frac\epsilon 2\,d(x,y)+\frac\epsilon 2\,d(z,t)\right) \\
& & \qquad\leq\;\exp\left(\frac\epsilon 2\,d(x,z)+\frac\epsilon 2\,d(y,t)\right)\:+\:
\exp\left(\frac\epsilon 2\,d(x,t)+\frac\epsilon 2\,d(y,z)\right)
\ee
for any $x,y,z,t\in X$.

\end{defn}
Note that if a space is Gromov hyperbolic with parameter $\delta$ then
the property also holds for any parameter $\wt\delta>\delta$.

In \cite{ZX} the following construction was made.  For a Ptolemaic space
$(X,d)$ and $U\subset X$ open with non-empty boundary let,
\begin{eqnarray}
\rho_U(x,y)= \sup_{p\in\partial U}\log
	\left(1+\frac{d(x,y)}{d(p,x)\,d(p,y)}\right).
\label{themetric}
\end{eqnarray}
This should be compared to \cite{AIW} where metrics hyperbolic in the
sense of Gromov were constructed on the complement of a
finite subset of a general metric space with parameter of hyperbolicity
independent of the size of the set.  Related metrics are to be found in
\cite{DHV}, \cite{PH}, and \cite{ZI} on domains in $\mathbb R^n$.  The
construction can also be seen as a variation of the inversion on metric
spaces (and the metric defined using it) in \cite{BHX}.

It was shown in \cite{ZX} that (\ref{themetric}) defines a metric on 
any open $U\subset X$ for any Ptolemaic space $(X,d)$.  Furthermore, the
following theorem was proved.

\begin{theorem}{\rm (\cite{ZX} Theorem 6)}

Let $(X,d)$ be a Ptolemaic space and $U\subset X$ open with non-empty
boundary.  If the distance between any two distinct points in $\partial U$
is at least $R>0$, then $(U,\rho_U)$ is Gromov hyperbolic with
parameter $\frac 12\log\max\left\{2+\frac{20}R,392\right\}$.
\label{ZXtheorem6}
\end{theorem}

When $(X,d)$ is Ptolemaic and $U$ is the complement of a point,
it was shown in \cite{ZX} (Theorem 5) that the resulting metric is strongly
hyperbolic with parameter $2$ and Gromov hyperbolic with parameter
$\frac 12\log 2$.  They also show (Theorem 4) that for $(X,d)$ Ptolemaic,
and $p\in X$ fixed, if
$$
s_p(x,y)=\frac{d(x,y)}{(1+d(x,p))(1+d(y,p))},
$$
then $(X,\log(1+s_p))$ is strongly hyperbolic with parameter $2$ and Gromov
hyperbolic with parameter $\frac 12\log 2$.

In this paper \thmref{ZXtheorem6} is generalized and sharpened.  It is shown
that the metric (\ref{themetric}) is Gromov hyperbolic with parameter $\log 2$
for any open $U$ whose boundary contains at least two points.  The possibility
is raised in \cite{ZX} that the construction may yield a strongly hyperbolic
metric.  It is proved (\thmref{maintheorem}) that this is indeed the case
with parameter $1$.  Examples are given (\thmref{hadamard}) showing that the
parameter bounds are sharp.

The paper is organized as follows.  In the first section it is proved
that for any open subset $U$ of a Ptolemaic space the construction produces
a metric that is strongly hyperbolic with parameter $1$, and Gromov hyperbolic
with parameter $\log 2$.  The second section is devoted to examples on
Hadamard manifolds exhibiting the sharp lower bound for the parameter of
Gromov hyperbolicity.

\section{Hyperbolic metrics}

In this section it is shown that the construction yields a strongly
hyperbolic metric with parameter $1$.

\begin{theorem}
If $(X,d)$ is a Ptolemaic space and $U\subset X$ is open with non-empty
boundary, then $(U,\rho_U)$ is Gromov hyperbolic with parameter $\log 2$,
and strongly hyperbolic with parameter $1$.  This result is sharp in the
sense that it does not hold in general for stronger assumptions on the
parameters when $\partial U$ contains at least two distinct points.
\label{maintheorem}
\end{theorem}

\noindent
Before proving this theorem a rearrangement inequality is proved as a lemma.

\begin{lemma}

If $\alpha,\beta,\gamma,\delta\in[0,\infty)$, then
\be
\min\{\alpha+\beta,\gamma+\delta\}\cdot\min\{\alpha+\gamma,\beta+\delta\} 
&\leq& \alpha\delta+\beta\gamma + 2\sqrt{\alpha\beta\gamma\delta}.
\ee
Equality is realized iff one or more of the following non-exclusive conditions
hold:
$$
\begin{array}{cl}
(\rc) & \alpha\delta=0\;\mbox{and}\;\max\{\alpha,\delta\}\geq|\beta-\gamma|,
  \\
(\rc) & \beta\gamma=0\;\mbox{and}\;\max\{\beta,\gamma\}\geq|\alpha-\delta|,
  \\
(\rc) & \alpha=\delta\;\mbox{and}\;\beta=\gamma.
\end{array}
$$
\label{thelemma}
\end{lemma}
\pf

\bigskip\noindent
Consider the case
\begin{eqnarray}
\alpha+\beta\leq\gamma+\delta\quad\mbox{and}\quad\alpha+\gamma\leq\beta+\delta,
\label{case1}
\end{eqnarray} 
Since the inequality as well as conditionis (i), (ii) and (iii) are invariant
under the permutations $(\alpha\:\delta)$, $(\alpha\:\beta)(\gamma\:\delta)$,
and $(\alpha\:\gamma\:\delta\:\beta)$, the other cases follow from this one.
From (\ref{case1}),
\begin{eqnarray}
\alpha &\leq& \delta-|\beta-\gamma| \label{case1a}
\end{eqnarray}
and the inequality becomes
\begin{eqnarray}
&\mbox{and}& (\alpha+\beta)(\alpha+\gamma)\;\leq\;\alpha\delta+\beta\gamma
	+2\sqrt{\alpha\beta\gamma\delta} \label{case1b} \\
&\Leftrightarrow& \sqrt\alpha\left(\sqrt\alpha^{\,3}
	+(\beta+\gamma-\delta)\sqrt\alpha
	-2\sqrt{\beta\gamma\delta}\right) \;\leq\; 0.
\label{case1c}
\end{eqnarray}
If $\delta=0$, then (\ref{case1a})$\;\Rightarrow
0\leq\alpha\leq-|\beta-\gamma|$, so $\beta=\gamma$ and
both sides of the inequality are equal to $\beta^2$.  If $\beta=0$, then
the (\ref{case1b}) becomes $\alpha(\alpha+\gamma)\leq\alpha\delta$, which is
true by (\ref{case1}).  The inequality holds similarly if $\gamma=0$.

Now assume $\beta\gamma\delta>0$.  Multiplying or dividing each side of
(\ref{case1b}) by $\delta^{\,2}$ does not change its validity, so without
loss of generality replace $(\alpha,\beta,\gamma,\delta)$ with
$(\alpha/\delta,\beta/\delta,\gamma/\delta,1),$ or simply assume
$\delta\geq 1$.  We claim that
$$
\phi(x)=x^3+(\beta+\gamma-\delta)x-2\sqrt{\beta\gamma\delta}
$$
is non-positive for $0\leq x<\sqrt{\delta-|\beta-\gamma|}$,  
from which (\ref{case1b}) follows.  Since $\phi(0)<0$, and $\phi$
has only one critical point in $[0,\infty)$, it's enough to show that
$\phi(\sqrt{\delta-|\beta-\gamma|})\leq 0$, that is
$$
\sqrt{\delta-|\beta-\gamma|}\,(\delta-|\beta-\gamma|+\beta+\gamma-\delta)
	\;\leq\; \sqrt{\delta}\,2\min\{\beta,\gamma\}
	\;\leq\; 2\sqrt{\beta\gamma\delta} 
$$
which is true when $\delta\geq 1$.

If equality holds in (\ref{case1b}), then
\be
\alpha=0 &\mbox{or}&
	\sqrt{\alpha}^{\,3}+(\beta+\gamma-\delta)\sqrt{\alpha}
	-2\sqrt{\beta\gamma\delta}=0.
\ee
If $\alpha=0$, then $\max\{\alpha,\delta\}=\delta\geq|\beta-\gamma|$, which
is covered by (i).  Otherwise, for reasons stated above, the only possibility
of a root of the cubic polynomial could be when $\alpha=\delta-|\beta-\gamma|$,
in which case
$$
\begin{array}{cccc}
 & \sqrt{\delta-|\beta-\gamma|}\,(\delta-|\beta-\gamma|+\beta+\gamma-\delta)
	&=& 2\sqrt{\beta\gamma\delta} \\
\Rightarrow& \sqrt{\delta-|\beta-\gamma|}\,\min\{\beta,\gamma\}
	&=& \sqrt{\beta\gamma\delta}.
\end{array}
$$
This equation is invariant under the permutation of $\beta$ with $\gamma$,
so without loss of generality assume that $\beta\leq\gamma$.  This gives
$$
(\delta-\gamma+\beta)\beta^2 \;=\; \beta\gamma\delta,
$$
so either $\beta=0$, in which case
$\alpha=\delta-\gamma=\delta-\max\{\beta,\gamma\}$, which is covered
by (ii), or $\beta>0$ and,
\be
 (\delta-\gamma+\beta)\beta = \gamma\delta 
&\Rightarrow& \beta = \gamma\;\Rightarrow\; \alpha = \delta,
\ee
which is (iii).

It is simple to verify that each of the conditions (i), (ii), and (iii)
yield equality.

\hfill\qed

\bigskip
Now to prove \thmref{maintheorem}.  Let $x,y,z,t\in U$ and $p,q\in\partial U$.
Since $(X,d)$ is Ptolemaic,
\be
d(x,y)\,d(z,p) &\leq&  d(x,z)\,d(y,p)+d(x,p)\,d(y,z), \\
d(x,y)\,d(t,p) &\leq&  d(x,t)\,d(y,p)+d(x,p)\,d(y,t), \\
d(z,t)\,d(x,q) &\leq&  d(x,z)\,d(t,q)+d(z,q)\,d(x,t), \\
\mbox{and}\quad d(z,t)\,d(y,q) &\leq&  d(y,z)\,d(t,q)+d(z,q)\,d(y,t).
\ee
Therefore,
\be
\frac{d(x,y)}{d(x,p)\,d(y,p)} &\leq & \frac{d(x,z)}{d(x,p)\,d(z,p)}+
	\frac{d(y,z)}{d(y,p)\,d(z,p)}, \\
\frac{d(x,y)}{d(x,p)\,d(y,p)} &\leq & \frac{d(x,t)}{d(x,p)\,d(t,p)}+
	\frac{d(y,t)}{d(y,p)\,d(t,p)}, \\
\frac{d(z,t)}{d(z,q)\,d(t,q)} &\leq & \frac{d(x,z)}{d(z,q)\,d(x,q)}+
	\frac{d(x,t)}{d(t,q)\,d(x,q)}, \\
\mbox{and}\quad\frac{d(z,t)}{d(z,q)\,d(t,q)} &\leq &
	\frac{d(y,z)}{d(z,q)\,d(y,q)}+ \frac{d(y,t)}{d(t,q)\,d(y,q)}.
\ee
For $a,b\in U$ denote the supremal metric space inversion over $\partial U$
by,
\begin{eqnarray}
\lambda(a,b)=\sup_{u\in\partial U}\frac{d(a,b)}{d(a,u)\,d(b,u)}.
\label{ld}
\end{eqnarray}
All the terms are positive so
\be
& & \mbox{$\left(1+\frac{d(x,y)}{d(x,p)\,d(y,p)}\right)\left(1+
	\frac{d(z,t)}{d(z,q)\,d(t,q)}\right)$} \\
&\leq& \mbox{$\min\left\{1+\frac{d(x,z)}{d(x,p)\,d(z,p)}
	+\frac{d(y,z)}{d(y,p)\,d(z,p)}\, ,\:
	1+\frac{d(x,t)}{d(x,p)\,d(t,p)}+\frac{d(y,t)}{d(y,p)\,d(t,p)}
	\right\}$} \\
& & \cdot\mbox{$\min\left\{1+\frac{d(x,z)}{d(z,q)\,d(x,q)}
	+\frac{d(x,t)}{d(t,q)\,d(x,q)}\, ,\:
	1+\frac{d(y,z)}{d(z,q)\,d(y,q)}+\frac{d(y,t)}{d(t,q)\,d(y,q)}
	\right\}$} \\
&\leq& \min\left\{1+\lambda(x,z)+\lambda(y,z)\, ,\:
	1+\lambda(x,t)+\lambda(y,t)\right\} \\
& & \cdot\min\left\{1+\lambda(x,z)+\lambda(x,t)\, ,\:
	1+\lambda(y,z)+\lambda(y,t)\right\} \\
&<& \min\left\{2+\lambda(x,z)+\lambda(y,z)\, ,\:
	2+\lambda(x,t)+\lambda(y,t)\right\} \\
& & \cdot\min\left\{2+\lambda(x,z)+\lambda(x,t)\, ,\:
	2+\lambda(y,z)+\lambda(y,t)\right\} .
\ee
Since $p,q\in\partial U$ were arbitrary,
\begin{eqnarray}
& & (1+\lambda(x,y))(1+\lambda(z,t)) \nonumber \\
&<& \min\left\{2+\lambda(x,z)+\lambda(y,z)\, ,\:
	2+\lambda(x,t)+\lambda(y,t)\right\} \nonumber \\
& & \cdot\min\left\{2+\lambda(x,z)+\lambda(x,t)\, ,\:
	2+\lambda(y,z)+\lambda(y,t)\right\} \nonumber \nonumber \\
&\leq& (1+\lambda(x,z))\,(1+\lambda(y,t))+(1+\lambda(y,z))\,(1+\lambda(x,t))
	\nonumber \\
& & +\; 2\sqrt{(1+\lambda(x,z))\,(1+\lambda(y,z))\,(1+\lambda(x,t))
	\,(1+\lambda(y,t))}.
\label{shest}
\end{eqnarray}
by \lemref{thelemma}.  From this it is easy to that,
\be
& & (1+\lambda(x,y))(1+\lambda(z,t)) \;< \\
& & \qquad 4\max\left\{(1+\lambda(x,z))(1+\lambda(y,t)),
	\:(1+\lambda(y,z))(1+\lambda(x,t))\right\}.
\ee
Gromov hyperbolicity of $(U,\rho_U)$ with parameter $\log 2$ follows
directly.  To show strong hyperbolicity, from (\ref{shest})
\be
\sqrt{(1+\lambda(x,y)(1+\lambda(z,t))} 
	&<& \sqrt{(1+\lambda(x,z))(1+\lambda(y,t))} \\
& & \quad +\;\sqrt{(1+\lambda(x,t)(1+\lambda(y,z))}.
\ee
Therefore $\rho_U$ is strongly hyperbolic with parameter $1$.

To see that the parameter is sharp, from \thmref{hadamard} it
follows that in general the parameter of Gromov hyperbolicity
can not be smaller than $\log 2$.  This and Theorem 4.2
in \cite{NS} give that the parameter of strong hyperbolicity
can not be larger than $1$.

\hfill\qed

\section{Examples on Hadamard manifolds}
\label{examples}

For a complete Riemannian manifold the Ptolemaic property is equivalent
to the manifold being Hadamard (Theorem 1.1 in \cite{BFW}).  This is
the setting for a class of examples which are Gromov hyperbolic with
parameter $\log 2$ and no lower.

\begin{theorem}
Let $(M,g)$ be a complete Riemannian manifold of dimension at least two
whose distance is Ptolemaic.  Take $U\subset M$ open with $p,q\in\partial U$
such that for some $R>0$
$$
d(p,q')\geq d(p,q)>0\:\mbox{ and }\: d(q,q')\geq R
$$
for all $q'\in\partial U\setminus\{p\}$.  Then $(U,\rho_U)$ is Gromov
hyperbolic with parameter $\log 2$ and not lower.
\label{hadamard}
\end{theorem}
\pf

Applying \thmref{maintheorem}, $(U,\rho_U)$ is $\log 2$-Gromov hyperbolic.
In the proof four points $x_+,x_-,y_+,y_-$ will be constructed with the
property that
\be
& & \rho_U(x_+,x_-)+\rho_U(y_+,y_-) \\
& & \qquad\leq\max\left\{\rho_U(x_+,y_+)+\rho_U(x_-,y_-),\:
  \rho_U(x_+,y_-)+\rho_U(x_-,y_+)\right\}+2\delta
\ee
with $\delta$ arbitrarily close to $\log 2$.

If $d_g(p,q)=2r$ take $\gamma\maps{[-r,r]}{M}$ a unit speed, minimal geodesic
with $\gamma(-r)=p$ and $\gamma(r)=q$.  The dimension of $M$ is at least two,
so we can construct a parallel unit vector field $F$ along
$\gamma$ with $\g{F}{\gamma'}=0$.  Fix $\theta>0$ small.

Define the curves
$$
\alpha_\theta(s)=\exp_{\gamma(-r\cos\theta)}(s F(-r\cos\theta))
\;\mbox{ and }\;
\beta_\theta(s)=\exp_{\gamma(r\cos\theta)}(s F(r\cos\theta))
$$
on the interval $[-r\sin\theta,r\sin\theta]$.  Let
$$
x_\pm=\alpha_\theta(\pm r\sin\theta),\; x_0=\alpha_\theta(0),
\; y_\pm=\beta_\theta(\pm r\sin\theta)\;\mbox{ and }\; y_0=\beta_\theta(0).
$$

\begin{figure}
\begin{center}
\includegraphics[scale=.9]{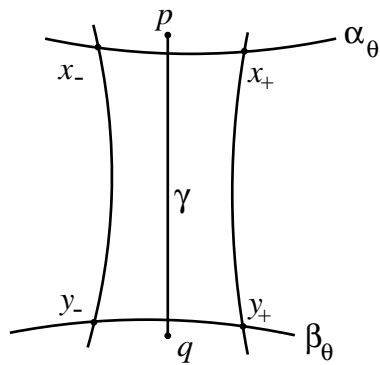}\\
\caption{A construction in a Hadamard manifold.}
\label{hfig}
\end{center}
\end{figure}

By construction,
$$
\g{\alpha_\theta'(0)}{\gamma'(-r\cos\theta)}=0 \qquad
	\g{\beta_\theta'(0)}{\gamma'(r\cos\theta)}=0.
$$
Since $(M,g)$ is complete, there exists $\kappa>0$ such that the absolute
value of the sectional curvature of $M$ is bounded above by $\kappa$ on a
closed Riemannian metric ball of radius $2r+R$ about $\gamma(0)$.  For all
$\theta>0$ sufficiently small, $r\sin\theta$ is less than the injectivity
radius at $\gamma(t)$ for all $-r\leq t\leq r$.  Take such a $\theta$.
Then there are unique minimal geodesics joining $x_0$, $p$,
and $x_\pm$ and joining $y_0$, $q$, and $y_\pm$ to form four geodesic
triangles.  Appling the Toponogov Comparison Theorem (Theorem 2.2 \cite{CE})
there are comparison triangles
$\triangle(\widetilde x_0,\widetilde p,\widetilde x_\pm)$, and
$\triangle(\widetilde y_0,\widetilde q,\widetilde y_\pm)$, in the
hyperbolic plane with constant curvature $-\kappa$ whose edges have the
same length as those in the corresponding triangle in $(M,g)$ and
whose corresponding angles are smaller.  Specifically,
\be
\angle(\widetilde p,\widetilde x_0,\widetilde x_\pm)
	\leq\angle(p,x_0,x_\pm)=\pi/2, & & 
\angle(\widetilde x_0,\widetilde p,\widetilde x_\pm)
	\leq\angle(x_0,p,x_\pm), \\
\angle(\widetilde q,\widetilde y_0,\widetilde y_\pm)
	\leq\angle(q,y_0,y_\pm)=\pi/2,&\mbox{and}&
\angle(\widetilde y_0,\widetilde q,\widetilde y_\pm)
	\leq\angle(y_0,q,y_\pm).
\ee
This and $\gamma|_{[-r\cos\theta,r]}$ and law of cosines in simply connected
space forms give (\ref{RII}).  
Applying the first and second variation formulas for length to variations
by geodesics along $\gamma|_{[-r\cos\theta,r\cos\theta]}$ gives
(\ref{sv1}); and along $\alpha_\theta|_{[-r\sin\theta,0]}$,
$\alpha_\theta|_{[0,r\sin\theta]}$, $\beta_\theta|_{[-r\sin\theta,0]}$, 
and $\beta_\theta|_{[0,r\sin\theta]}$ give (\ref{sv2}) where
$\sigma, \tau>0$ are independent of $\theta\in(0,\theta_0)$ for some
$\theta_0>0$.
\begin{eqnarray}
2r-\sigma\theta \;\leq & d_g(x_\pm,q),\, d_g(y_\pm,p) & \leq\;
	2r+\sigma\theta 
\label{RII} \\
2r\cos\theta-\tau\theta^2 \;\leq & d_g(x_\pm,y_\pm), d_g(x_\pm,y_\mp) &
	\leq\; 2r\cos\theta+\tau\theta^2 
	\label{sv1} \\
r\sin\theta-\tau\theta^2 \;\leq & d_g(x_\pm,p),\, d_g(y_\pm,q) & \leq\;
	r\sin\theta+\tau\theta^2
	\label{sv2} 
\end{eqnarray}
Also, by construction,
\begin{eqnarray}
d_g(x_-,x_+)=d_g(y_-,y_+)=2r\sin\theta.
\label{fc}
\end{eqnarray}
We claim that for $\theta>0$ sufficiently small, the suprema over
$\partial U$ in the definition (\ref{ld}) of
$\lambda(x_\pm,y_\pm)$, $\lambda(x_\pm,y_\mp)$, $\lambda(x_-,x_+)$, and
$\lambda(y_-,y_+)$ are realized at $p$ or at $q$.  This together with
(\ref{sv1}), (\ref{sv2}), (\ref{RII}), and (\ref{fc}) imply that
$$
\lim_{\theta\rightarrow 0^+}\mbox{$\frac
	{4\max\{(1+\lambda(x_-,y_-)(1+\lambda(x_+,y_+)),\:
	(1+\lambda(x_-,y_+)(1+\lambda(x_+,y_-))\}}
	{(1+\lambda(x_-,x_+))(1+\lambda(y_-,y_+))}
	=1$}
$$
from which it follows that $\rho_U$ is Gromov hyperbolic with parameter
$\log 2$ and no smaller.

Now to prove the claim.  Let
$$
\eta=\max\{d_g(x_\pm,p),d_g(y_\pm,q)\}<\min\{r/2,R/3\}.
$$
Then for $w\in\partial U\setminus\{p,q\}$,
\be
d_g(y_-,w) &\geq& d_g(q,w)-d_g(y_-,q)\;\geq\; R-\eta \;>\;\frac 23 R \\
d_g(x_-,w) &\geq& d_g(p,w)-d_g(x_-,p)\;\geq\;
	d_g(q,p)-\eta = 2r-\eta > 3\eta.
\ee
Therefore,
\be
d_g(x_-,p)\,d_g(y_-,p) &\leq& \eta\left((d_g(y_-,q)+d_g(q,p)\right)
	\;\leq\; \eta\left(\eta+d_g(p,w)\right) \\
&\leq& \eta\left(\eta+d_g(p,x_-)+d_g(x_-,w)\right)
	\;\leq\; \eta\left(2\eta+d_g(x_-,w)\right) \\
&<& \frac{5R}{9}\,d_g(x_-,w)
	\;<\;\frac 56 \, d_g(y_-,w)\,d_g(x_-,w).
\ee
The inequalities
\be
d_g(x_-,p)\,d_g(y_+,p) &<& \frac 56\,d_g(x_-,w)\,d_g(y_+,w), \\
d_g(x_+,p)\,d_g(y_-,p) &<& \frac 56\,d_g(x_+,w)\,d_g(y_-,w), \\
\mbox{and}\quad d_g(x_+,p)\,d_g(y_+,p) &<& \frac 56\,d_g(x_+,w)\,d_g(y_+,w), 
\ee
are proved similarly.  From this we obtain
\be
d_g(x_-,q)\,d_g(y_-,q) &<& d_g(x_-,w)\,d_g(y_-,w), \\
d_g(x_-,q)\,d_g(y_+,q) &<& d_g(x_-,w)\,d_g(y_+,w), \\
d_g(x_+,q)\,d_g(y_-,q) &<& d_g(x_+,w)\,d_g(y_-,w), \\
\mbox{and}\quad d_g(x_+,q)\,d_g(y_+,q) &<& d_g(x_+,w)\,d_g(y_+,w),
\ee
using (\ref{sv2}) and (\ref{RII}) for all $w\in\partial U$ when $\theta>0$
is sufficiently small.  For the other cases,
\be
d_g(x_+,w) &\geq& d_g(p,w)-d_g(x_+,p)\;\geq\; 2r-\eta \;>\; 3\eta, \\
\mbox{and}\quad d_g(x_-,w) &\geq& d_g(p,w)-d_g(x_-,p)\;\geq\; 2r-\eta \;>\;
	3\eta.
\ee
Therefore,
$$
d_g(x_-,p)\,d_g(x_+,p) \;\leq\; \eta^2 < \frac19\,d_g(x_-,w)\,d_g(x_+,w).
$$
The inequality
$$
d_g(y_-,p)\,d_g(y_+,p) < \frac19\,d_g(y_-,w)\,d_g(y_+,w),
$$
is obtained in a similar way.  Therefore for $z_1,z_2\in\{x_\pm,y_\pm\}$
distinct we have that
$$
\max\left\{\frac{d_g(z_1,z_2)}{d_g(z_1,p)\,d_g(z_2,p)},
	\frac{d_g(z_1,z_2)}{d_g(z_1,q)\,d_g(z_2,q)}\right\}
	\geq \frac{d_g(z_1,z_2)}{d_g(z_1,w)\,d_g(z_2,w)}
$$
for all $w\in\partial U$ which was the claim.

\hfill\qed

\bigskip
\noindent
Mathematics Department, New York City College of Technology, 300 Jay Street,
Brooklyn NY, USA 11201.

\end{document}